\documentclass[a4paper,11pt]{amsart}

\usepackage[margin=0.95in]{geometry}
\usepackage{amsmath,amssymb,mathtools}
\usepackage{booktabs,tabularx,array,longtable}
\usepackage{enumitem}
\usepackage{microtype}
\usepackage{needspace}
\usepackage{xcolor}
\usepackage[colorlinks=true,linkcolor=blue!45!black,citecolor=blue!45!black,
  urlcolor=blue!45!black]{hyperref}
\usepackage[capitalize,noabbrev]{cleveref}
\newtheorem{theorem}{Theorem}
\newtheorem{proposition}[theorem]{Proposition}
\newtheorem{corollary}[theorem]{Corollary}
\newtheorem{lemma}[theorem]{Lemma}
\theoremstyle{definition}

\newtheorem{example}[theorem]{Example}
\theoremstyle{remark}

\crefname{theorem}{theorem}{theorems}
\Crefname{theorem}{Theorem}{Theorems}
\crefname{proposition}{proposition}{propositions}
\Crefname{proposition}{Proposition}{Propositions}
\crefname{corollary}{corollary}{corollaries}
\Crefname{corollary}{Corollary}{Corollaries}
\crefname{example}{example}{examples}
\Crefname{example}{Example}{Examples}

\AddToHook{env/theorem/begin}{\crefalias{theorem}{theorem}}
\AddToHook{env/proposition/begin}{\crefalias{theorem}{proposition}}
\AddToHook{env/corollary/begin}{\crefalias{theorem}{corollary}}
\AddToHook{env/lemma/begin}{\crefalias{theorem}{lemma}}
\AddToHook{env/definition/begin}{\crefalias{theorem}{definition}}
\AddToHook{env/example/begin}{\crefalias{theorem}{example}}
\AddToHook{env/remark/begin}{\crefalias{theorem}{remark}}

\newcommand{\CC}{\mathbb C}
\newcommand{\FF}{\mathbb F}
\newcommand{\QQ}{\mathbb Q}
\newcommand{\RR}{\mathbb R}
\newcommand{\ZZ}{\mathbb Z}
\newcommand{\PP}{\mathbb P}
\newcommand{\Fix}{\operatorname{Fix}}
\newcommand{\Res}{\operatorname{Res}}
\newcommand{\rk}{\operatorname{rk}}
\newcommand{\tr}{\operatorname{tr}}
\newcommand{\eul}{\operatorname e}

\title[Fixed-locus relations on $K3$ surfaces]
{Fixed-locus relations and purely non-symplectic automorphisms of order $6$ on
$K3$ surfaces}
\author{Dominik Burek}
\address{Institute of Mathematics, Jagiellonian University, ul. \L{}ojasiewicza 6,
30-348 Krak\'ow, Poland}
\email{dominik.burek@uj.edu.pl}
\hypersetup{
  pdftitle={Fixed-locus relations and purely non-symplectic automorphisms of order 6 on K3 surfaces},
  pdfauthor={Dominik Burek},
  pdfsubject={Purely non-symplectic automorphisms of K3 surfaces},
  pdfkeywords={K3 surface, non-symplectic automorphism, fixed locus, Borcea--Voisin construction}
}

\begin{document}

\begin{abstract}
We study relations between invariants of the fixed loci of purely
non-symplectic automorphisms of $K3$ surfaces.  We obtain these relations by
comparing the stringy Euler characteristic of higher dimensional Borcea--Voisin
varieties with their Hodge numbers.  For orders $4$ and $6$ we explain the
geometric meaning of these relations and show how they follow from the
Lefschetz and Riemann--Hurwitz formulas.  For order $6$, the classifications of
fixed loci for orders $2$ and $3$ give $817$ candidates.  We use local conditions
on discriminant forms, Hermitian trace forms and bounds for root-free lattices
to reduce this number to $204$.  All $142$ numerical types obtained from the
$150$ deformation classes of Brandhorst and Hofmann satisfy these conditions.
Also, we give explicit elliptic and projective models for each of the $20$
possible fixed loci of the generator.
\end{abstract}

\maketitle

\section{Introduction}

One of the reasons for studying non-symplectic automorphisms of $K3$ surfaces is
the construction of Calabi--Yau manifolds due to Borcea and
Voisin~\cite{Borcea,Voisin}.  They used a non-symplectic involution of a $K3$
surface and the involution of an elliptic curve.  Cattaneo and Garbagnati
extended this construction to automorphisms of orders $3$, $4$ and $6$.
Combining their method with the construction of Cynk and Hulek gives a higher
dimensional construction
\cite{CattaneoGarbagnati,CynkHulek,BurekBV}.

Let $S$ be a complex $K3$ surface and let $\sigma\in\operatorname{Aut}(S)$ have
finite order $d$.  The automorphism $\sigma$ is called \emph{purely
non-symplectic} if
\[
 \sigma^*\omega_S=\zeta_d\omega_S,
\]
where $\zeta_d$ is a primitive $d$th root of unity and
$0\ne\omega_S\in H^{2,0}(S)$.  The fixed loci of $\sigma$ and its powers are
unions of smooth curves and isolated points.  Their genera, orbit lengths and local tangent types enter the formulas for
the Hodge numbers of the resulting Calabi--Yau manifolds.

For $d\in\{2,3,4,6\}$, let $E_d$ be an elliptic curve with an automorphism $\tau$ of
order $d$ satisfying $\tau^*\omega_{E_d}=\zeta_d\omega_{E_d}$, and put
\[
G_{d,n}=\{(a_0,\ldots,a_{n-1})\in(\ZZ/d\ZZ)^n:
             a_0+\cdots+a_{n-1}=0\}.
\]
The group acts on the factors through $\sigma$ and $\tau$, respectively.
The quotient of $S\times E_d^{n-1}$ by $G_{d,n}$ admits a crepant resolution
$Y_{d,n}$ of dimension $n+1$ under the hypotheses of~\cite{BurekBV}.  We can compute its Euler characteristic
using either the orbifold formula or its Hodge polynomial.  If a crepant
resolution exists, the two computations give the same result~\cite{DixonEtAl1,DixonEtAl2,BatyrevStringy}.

The complete classification of finite groups acting on complex $K3$ surfaces
is due to Brandhorst and Hofmann~\cite{BrandhorstHofmann,BrandhorstHofmannData}.
Their list contains $150$ deformation classes of purely non-symplectic
automorphisms of order $6$.  Our contribution concerns explicit equations for
the fixed loci and the information obtained from the stringy Euler calculations.

We give equations for all $20$ generator loci, combining Dillies' and
Artebani--Sarti's constructions with further elliptic and projective
models~\cite{Dillies6,ArtebaniSarti3}.
Example~\ref{ex:missing-six} gives the model $y^2=x^3+t^6x+1$ with
$\sigma(x,y,t)=(x,y,\zeta_6t)$ for the elliptic-curve case with four isolated
points identified by Brandhorst and Hofmann and missing from Dillies' list.
We compute the action on the resolution and the fixed loci of both proper powers.

The comparison with Hodge numbers gives relations between fixed-locus
invariants.  For order $4$,
Theorem~\ref{thm:order4} identifies their quotient-genus terms.  For order $6$,
Propositions~\ref{prop:balance} and \ref{prop:dependence} express these relations
in terms of Lefschetz, Riemann--Hurwitz, the dimension formula and the orbit
decomposition.

Using the order-$2$ and order-$3$ classifications, we reduce $817$ candidates
to $204$ by a finite lattice sieve.  We prove Euler recurrences in all dimensions.
For order $6$, the full sequence determines, and is determined by, the three
fixed-locus Euler characteristics (Proposition~\ref{prop:euler-determination}).
This gives $102$ Euler sequences in the Brandhorst--Hofmann list.
We obtain the following results.

\Needspace{8\baselineskip}
\begin{theorem}\label{thm:introduction}
Let $\sigma$ be a purely non-symplectic automorphism of a complex $K3$ surface.
\begin{enumerate}[label=\textup{(\roman*)}]
\item Suppose that $\sigma$ has order $4$ and $\Fix(\sigma^2)$ is not the
union of two elliptic curves.  If a $\sigma$-invariant maximal curve is not
fixed pointwise, two dependent linear relations determine its quotient genus.
\item If $\sigma$ has order $6$, the stringy relations considered here follow
from the dimension, Lefschetz, orbit and Riemann--Hurwitz relations.  We give
explicit identities expressing these dependencies.
\item The $150$ deformation classes of Brandhorst and Hofmann give $142$
numerical fixed-locus types, $20$ possible fixed loci of $\sigma$, and $102$
different sequences of stringy Euler characteristics.  Each of the $20$ fixed
loci of $\sigma$ is realized below by an explicit equation.
\item Starting with the order-$2$ and order-$3$ fixed loci gives $817$
possibilities.  Necessary lattice conditions reduce this number to $204$.
Every one of the $142$ realized numerical types belongs to this shorter list.
\end{enumerate}
\end{theorem}

The sieve does not provide an equivariant primitive embedding in the $K3$
lattice.  We use the theorem of Brandhorst and Hofmann for completeness.

In Section~2 we recall the Lefschetz formulas and in Section~3 we compute the
stringy Euler characteristic.  Sections~4 and~5 concern automorphisms of orders
$4$ and $6$.  The explicit constructions are given at the end of Section~5.
In Section~6 we study the two proper powers of an automorphism of order $6$.
The last section contains the application to the complete list and to
Borcea--Voisin manifolds.

\section{Fixed-locus and cohomological data}

For a finite-order automorphism $g$ of a $K3$ surface, the topological Lefschetz
formula gives
\begin{equation}\label{eq:top-lef}
  2+\tr\bigl(g^*\mid H^2(S,\CC)\bigr)=\eul(\Fix(g)).
\end{equation}
If the fixed locus consists of curves $C_j$ and isolated points $P_i$, then
\(\eul(\Fix(g))=\sum_j(2-2g(C_j))+\#\{P_i\}\).
We use the holomorphic Lefschetz formula only for the structure sheaf.  Its local
terms at a fixed point with tangent eigenvalues $(\lambda_1,\lambda_2)$ and at a
fixed curve $C$ with normal eigenvalue $\lambda$ are
\begin{equation}\label{eq:hol-local}
 \frac{1}{(1-\lambda_1)(1-\lambda_2)},\qquad
 \frac{1-g(C)}{1-\lambda}-\frac{\lambda C^2}{(1-\lambda)^2},
\end{equation}
respectively.  Adjunction gives $C^2=2g(C)-2$.

If $C\subset\Fix(g^j)$, the generator $g$ may fix $C$ pointwise, act
non-trivially on $C$, or permute it with other components.  In the second case
we use the quotient genus, and in the third case we count the orbits.

\section{The stringy Euler characteristic}

Let a finite group $G$ act algebraically on a smooth projective complex
variety $X$.  The physicists' stringy, or orbifold,
Euler number is
\begin{equation}\label{eq:orbifold-euler}
 \eul_{\mathrm{orb}}(X/G)
 =\frac{1}{|G|}\sum_{\substack{g,h\in G\\gh=hg}}
       \eul(X^g\cap X^h).
\end{equation}
If the action preserves the volume form and $X/G$ admits a crepant resolution
$\widetilde{X/G}$, then
\begin{equation}\label{eq:stringy-crepant}
 \eul(\widetilde{X/G})=\eul_{\mathrm{orb}}(X/G)
\end{equation}
in the present abelian setting
\cite{DixonEtAl1,DixonEtAl2,BatyrevStringy}.  Equation
\eqref{eq:orbifold-euler} can be applied directly, since every common fixed locus is a
product of fixed loci on $S$ and $E_d$.

Put
\[
 \varepsilon_j=\eul(\Fix(\sigma^j)),\qquad \varepsilon_d=\eul(S)=24.
\]
For $x=(a,b)\in A_d:=(\ZZ/d\ZZ)^2$, set
\(\nu(x)=\gcd(d,a,b)\).  If $\tau$ is the standard automorphism of $E_d$,
write
\[
 u_d(x)=\eul\bigl(E_d^{\tau^a}\cap E_d^{\tau^b}\bigr).
\]
The only values needed are
\begin{center}
\begin{tabular}{@{}c|ccc@{}}
\toprule
$d$ & \multicolumn{3}{c}{$u_d(x)$ according to $\nu(x)$}\\
\midrule
$4$ & $u_4(1)=2$ & $u_4(2)=4$ & $u_4(4)=0$\\
$6$ & $u_6(1)=1$ & $u_6(2)=3$ & $u_6(3)=4$, $u_6(6)=0$\\
\bottomrule
\end{tabular}
\end{center}

\begin{theorem}\label{thm:stringy-convolution}
For $d\in\{4,6\}$ and every $n\geq1$,
\begin{equation}\label{eq:stringy-convolution}
 \eul(Y_{d,n})=
 \frac{1}{d^{n-1}}
 \sum_{\substack{x_0,\ldots,x_{n-1}\in A_d\\
                  x_0+\cdots+x_{n-1}=0}}
 \varepsilon_{\nu(x_0)}\prod_{j=1}^{n-1}u_d(x_j).
\end{equation}
For the first six values we obtain the following formulas.
\begin{center}
\small
\begin{tabularx}{\textwidth}{@{}cXX@{}}
\toprule
$n$ & $\eul(Y_{4,n})$ & $\eul(Y_{6,n})$\\
\midrule
$1$ & $\varepsilon_4$
    & $\varepsilon_6$\\
$2$ & $6\varepsilon_1+3\varepsilon_2$
    & $4\varepsilon_1+4\varepsilon_2+2\varepsilon_3$\\
$3$ & $60\varepsilon_1+15\varepsilon_2+6\varepsilon_4$
    & $64\varepsilon_1+24\varepsilon_2+8\varepsilon_3+4\varepsilon_6$\\
$4$ & $546\varepsilon_1+138\varepsilon_2+45\varepsilon_4$
    & $660\varepsilon_1+228\varepsilon_2+84\varepsilon_3+28\varepsilon_6$\\
$5$ & $4920\varepsilon_1+1230\varepsilon_2+411\varepsilon_4$
    & $6656\varepsilon_1+2232\varepsilon_2+832\varepsilon_3+280\varepsilon_6$\\
$6$ & $44286\varepsilon_1+11073\varepsilon_2+3690\varepsilon_4$
    & $66644\varepsilon_1+22244\varepsilon_2+8332\varepsilon_3+2780\varepsilon_6$\\
\bottomrule
\end{tabularx}
\end{center}
\end{theorem}

\begin{proof}
A commuting pair in the abelian group $G_{d,n}$ corresponds to a sequence
$x_0,\ldots,x_{n-1}\in A_d$ whose sum is zero.  The subgroup generated by
$\sigma^a$ and $\sigma^b$ is generated by $\sigma^{\gcd(d,a,b)}$; hence the
$K3$ factor contributes $\varepsilon_{\nu(x_0)}$.  The remaining factors
contribute $u_d(x_j)$.  Dividing their product sum by
$|G_{d,n}|=d^{n-1}$ proves \cref{eq:stringy-convolution}.  Computing the convolution of
$u_d$ on the finite group $A_d$ gives the coefficients in the table.
\end{proof}

We apply this formula to orders four and six.  For order
$4$, let $\Gamma_j$ be the sum of the genera of the fixed curves of $\sigma^j$.
Then
\[
 \varepsilon_1=s+2k-2\Gamma_1,\qquad
 \varepsilon_2=2N-2\Gamma_2,
\]
and therefore
\begin{equation}\label{eq:stringy4-threefold}
 \eul(Y_{4,2})=6s+12k-12\Gamma_1+6N-6\Gamma_2.
\end{equation}
If the maximal curve $D\subset\Fix(\sigma^2)$ is fixed pointwise, then
$\Gamma_1=\Gamma_2=g(D)$; if it moves, then its genus occurs only in
$\Gamma_2$.  We therefore obtain the two cases in
\cref{thm:order4}.  Substituting the corresponding Hodge formula, we find that
the Hodge Euler number minus \cref{eq:stringy4-threefold} is twice the left-hand
side of \cref{eq:R4b}.

For order $6$, write $\Gamma_j$ for the sum of the genera of all curves in
$\Fix(\sigma^j)$.  In the notation of \cref{sec:order6},
\begin{align*}
 \varepsilon_1&=p_{(2,5)}+p_{(3,4)}+2\ell-2\Gamma_1,\\
 \varepsilon_2&=n+2k-2\Gamma_2,\\
 \varepsilon_3&=2N-2\Gamma_3.
\end{align*}
Using $n=p_{(2,5)}+2n'$ gives the threefold calculation
\begin{align}\label{eq:stringy6-threefold}
 \eul(Y_{6,2})={}&8\ell-8\Gamma_1+8p_{(2,5)}+4p_{(3,4)}
 +8k-8\Gamma_2+8n'\nonumber\\
 &+4N-4\Gamma_3.
\end{align}
The difference between the Hodge Euler formula from
\cite[Proposition~7.3]{CattaneoGarbagnati} and
\cref{eq:stringy6-threefold} is $2\mathcal B_6$.  Thus we obtain \cref{eq:B6}.
In Proposition~\ref{prop:balance} we show that it follows from Lefschetz and
Riemann--Hurwitz for the residual actions.

\section{Relations for automorphisms of order four}

Let $\sigma^*\omega_S=i\omega_S$.  On cohomology write
\begin{equation}\label{eq:coh4}
H^2(S,\CC)=V_1^{\oplus r}\oplus V_i^{\oplus m}
 \oplus V_{-1}^{\oplus(22-r-2m)}\oplus V_{-i}^{\oplus m}.
\end{equation}
Assume that $\Fix(\sigma^2)$ is not the union of two elliptic curves.
Choose a $\sigma$-invariant component $D$ of maximal genus, put $g=g(D)$,
and let $N$ be the number of components.  Such a component exists: holomorphic
Lefschetz implies $\Fix(\sigma)\ne\varnothing$, so at least one component
of $\Fix(\sigma^2)$ is invariant.  If the maximal genus is positive, that
component is unique.  Following~\cite{CattaneoGarbagnati}, we divide these components
into three types:
\begin{center}
\begin{tabularx}{0.92\textwidth}{@{}cXc@{}}
\toprule
type & action of $\sigma$ on a component of $\Fix(\sigma^2)$ & number \\
\midrule
I & fixed pointwise & $k$ \\
II & invariant, but not fixed pointwise & $b$ \\
III & exchanged in pairs & $2a$ \\
\bottomrule
\end{tabularx}
\end{center}
Thus
\begin{equation}\label{eq:decomp4}N=k+b+2a.
\end{equation}
Let $n_1$ and $n_2$ count the isolated fixed points of $\sigma$ off and on $D$,
respectively, and put $s=n_1+n_2$.  Finally set
\[
h=\sum_{C\subset\Fix(\sigma)}(1-g(C)).
\]
The fixed-locus formulas of Artebani--Sarti~\cite{ArtebaniSarti4}, in the
notation used here, are
\begin{equation}\label{eq:AS4}
 r=\frac{12+k+2a+b-g+4h}{2},\qquad
 m=\frac{12-k-2a-b+g}{2},\qquad s=2h+4.
\end{equation}
In particular, topological Lefschetz for $\sigma^2$ is
\begin{equation}\label{eq:TL4square}N-g-12+2m=0.
\end{equation}

There are two cases.  If $D$ is of type I, then
\begin{equation}\label{eq:caseI}
h=k-g,\qquad n_2=0,\qquad b=\frac{n_1}{2}.
\end{equation}
If $D$ is of type II, the restriction of $\sigma$ to $D$ is an involution
branched at the $n_2$ fixed points, and
\begin{equation}\label{eq:caseII}
h=k,\qquad b=\frac{n_1}{2}+1,\qquad
q:=g(D/\sigma)=\frac{2g+2-n_2}{4}.
\end{equation}

Define $\delta=0$ in the type-I case and $\delta=q$ in the type-II case.
Put
\begin{align*}
R_{4,1}&=-b+8+3k-3g+2s+2a-2r,\\
R_{4,2}&=-m-2k+2g-s-2a+r+2.
\end{align*}

\begin{theorem}\label{thm:order4}
With the notation above,
\begin{align}
0={}&R_{4,1}+4\delta, \label{eq:R4a}\\
0={}&R_{4,2}-2\delta. \label{eq:R4b}
\end{align}
Moreover, these are dependent: the left-hand side of
\cref{eq:R4a} plus twice the left-hand side of \cref{eq:R4b} equals
the negative of the left-hand side of \cref{eq:TL4square}.
\end{theorem}

\begin{proof}
If $D$ is of type I, substitute \cref{eq:caseI} and $s=2h+4$ into
\cref{eq:AS4}; both expressions vanish.  Suppose that $D$ is of type II.
By Riemann--Hurwitz we have $n_2=2g+2-4q$, so
\[
s=2k+4,\qquad b=k+2-g+2q.
\]
Substitution into \cref{eq:AS4} gives $R_{4,1}=-4q$ and
$R_{4,2}=2q$.  This proves
\cref{eq:R4a,eq:R4b}.  The dependence follows directly from
\[
(-b+8+3k-3g+2s+2a-2r)
+2(-m-2k+2g-s-2a+r+2)=12-N+g-2m.
\]
\end{proof}

\medskip
\begin{corollary}\label{cor:qrecover}
If $D$ is of type II, either linear form determines its quotient genus:
\[
g(D/\sigma)=-\frac14R_{4,1}=\frac12R_{4,2}.
\]
Consequently $R_{4,1}=R_{4,2}=0$ precisely when $D$ is of type I or when $D$
is of type II with rational quotient.
\end{corollary}

\begin{example}\label{ex:counter4}
The first row of Artebani--Sarti's Table~5~\cite{ArtebaniSarti4} is realized by an order-$4$ action
with
\[
(r,m,k,g,n_1,n_2,a,b,N)=(2,10,0,10,2,2,0,2,2).
\]
Here $D$ is of type II and $g(D/\sigma)=5$.  We have
$R_{4,1}=-20$ and $R_{4,2}=10$, and substituting the quotient genus in \cref{eq:R4a,eq:R4b}
gives zero.
\end{example}

\subsection{Explicit \texorpdfstring{order-$4$}{order-4} models}

\begin{example}\label{ex:quartic}
Let $Q_4$ be a smooth plane quartic and consider
\[
S_Q:\quad w^4=Q_4(x_0,x_1,x_2),\qquad
\sigma(w,x_0,x_1,x_2)=(iw,x_0,x_1,x_2).
\]
This is the cyclic degree-$4$ cover of $\PP^2$ branched along $Q_4$.
The canonical bundle formula gives
\(K_{S_Q}=\pi^*(K_{\PP^2}+3\mathcal O_{\PP^2}(1))\simeq\mathcal O_{S_Q}\).
Since $S_Q$ is a smooth quartic in $\PP^3$, it also has
$H^1(S_Q,\mathcal O_{S_Q})=0$, and hence is a $K3$ surface.  The ramification curve is fixed
pointwise and has genus $3$.  From this model and \cref{eq:AS4} we get
\[
(r,m,k,N,a,b,g,n_1,n_2)=(1,7,1,1,0,0,3,0,0).
\]
The associated order-$4$ Borcea--Voisin threefold has Hodge numbers $(9,27)$.
\end{example}

\medskip
\begin{example}\label{ex:ell4}
Let $a_8(t)$ be a degree-$8$ polynomial with simple roots and take the elliptic
$K3$ surface
\[
S_a:\quad y^2=x^3-a_8(t)x,\qquad
\sigma(x,y,t)=(-x,iy,t).
\]
The curve $D:y=x^2-a_8(t)=0$ is a genus-$3$ component of
$\Fix(\sigma^2)$ on which $\sigma$ acts non-trivially.  The fixed-locus data are
\[
(r,m,k,N,a,b,g,n_1,n_2)=(10,6,2,3,0,1,3,0,8),
\]
so $q=(2g+2-n_2)/4=0$.  This example has a type-II maximal curve and
$R_{4,1}=R_{4,2}=0$; its associated threefold has Hodge numbers $(44,8)$.
Degenerating $a_8$ changes the ramification and gives
several rows of the order-$4$ tables~\cite{ArtebaniSarti4}.
\end{example}

The Hodge numbers used in these examples are
\cite[Section~6]{CattaneoGarbagnati}
\begin{align}\label{eq:hodge4}
h^{1,1}&=1+r+7k+3b+2(n_1+n_2)+4a,\nonumber\\
h^{2,1}&=
\begin{cases}
m-1+7g,&D\text{ of type I},\\
m+2g-n_2/2,&D\text{ of type II}.
\end{cases}
\end{align}
By Riemann--Hurwitz, the contribution of the moving curve to $h^{2,1}$ also
determines its quotient genus.

\section{Relations for automorphisms of order six}\label{sec:order6}

Let $\zeta=\zeta_6$ and suppose $\sigma^*\omega_S=\zeta\omega_S$.  Use the
eigenspace decomposition
\begin{equation}\label{eq:coh6}
H^2(S,\CC)=V_1^{\oplus r}\oplus
V_{\zeta,\zeta^5}^{\oplus m}\oplus
V_{\zeta^2,\zeta^4}^{\oplus\alpha}\oplus V_{-1}^{\oplus\beta},
\end{equation}
where each member of a conjugate pair has the indicated multiplicity.  Hence
\begin{equation}\label{eq:D6}
\mathcal D_6:=r+2m+2\alpha+\beta-22=0.
\end{equation}

We use the following notation for the fixed loci and the action of $\sigma$.
The locus $\Fix(\sigma)$ contains $\ell$ curves, one of maximal genus $D$, and
isolated points of tangent types $(\zeta^2,\zeta^5)$ and
$(\zeta^3,\zeta^4)$, counted by $p_{(2,5)}$ and $p_{(3,4)}$.  The locus
$\Fix(\sigma^2)$ contains $k$ curves, one of maximal genus $G$, and $n$ isolated
points.  Among its curves, $b$ pairs are exchanged by $\sigma$; among its isolated
points, $n'$ pairs are exchanged.  The remaining isolated points are the
$p_{(2,5)}$ points already fixed by $\sigma$, so
\begin{equation}\label{eq:P6}
\mathcal P_6:=n-p_{(2,5)}-2n'=0.
\end{equation}
Finally, $\Fix(\sigma^3)$ contains $N$ curves; $a$ triples are cyclically
permuted by $\sigma$, and $F_1,F_2$ denote the components of largest genera.
We set $g(D)=0$ if $\Fix(\sigma)$ has no curves and $g(G)=0$ if
$\Fix(\sigma^2)$ has no curves.  Missing components $F_j$ also contribute zero.
Thus these symbols record the sums of the genera in the three fixed loci.
The positive-genus components of $\Fix(\sigma^2)$ and $\Fix(\sigma^3)$ are
$\sigma$-invariant.  We write $q_G$ and $q_1+q_2$ for the sums of the genera
of their quotients; pointwise fixed components contribute their original genera.

Topological Lefschetz for the three non-trivial powers gives
\begin{align}
\mathcal T_1:={}&2+r+m-\alpha-\beta-2\ell+2g(D)
 -p_{(2,5)}-p_{(3,4)}=0,\label{eq:T1}\\
\mathcal T_2:={}&2+r-m-\alpha+\beta-n-2k+2g(G)=0,\label{eq:T2}\\
\mathcal T_3:={}&2+r+2\alpha-\beta-2m-2N
 +2g(F_1)+2g(F_2)=0.\label{eq:T3}
\end{align}
Holomorphic Lefschetz for
$\mathcal O_S$ gives
\begin{equation}\label{eq:H6}
3+3\ell-3g(D)-\frac12p_{(3,4)}-p_{(2,5)}=0.
\end{equation}
Riemann--Hurwitz, applied to the invariant components on which $\sigma$ acts
non-trivially, gives
\begin{align}
q_G&=\frac{2g(G)+2g(D)-p_{(3,4)}+2k-4b-2\ell}{4},\label{eq:RHG}\\
q_1+q_2&=\frac{2g(F_1)+2g(F_2)+4g(D)-2p_{(2,5)}-2p_{(3,4)}
 +4N-12a-4\ell}{6}.\label{eq:RHF}
\end{align}
To prove these formulas, we remove the $\ell$ pointwise fixed curves and the
permuted components.  For $\Fix(\sigma^2)$, there remain $k-2b-\ell$
curves of total genus $g(G)-g(D)$, with quotient genera summing to
$q_G-g(D)$.  The involution has $p_{(3,4)}$ fixed points on this union.
For $\Fix(\sigma^3)$, there remain $N-3a-\ell$ curves of total genus
$g(F_1)+g(F_2)-g(D)$, with quotient genera summing to $q_1+q_2-g(D)$.
The order-$3$ action has $p_{(2,5)}+p_{(3,4)}$ fixed points there, each with
ramification index $3$.  Summing Riemann--Hurwitz over these components
proves both formulas, including the cases $g(D)>0$.

Consider the following linear form for order $6$:
\begin{align}
\mathcal B_6:={}&-m+2+r-2\ell-p_{(2,5)}-p_{(3,4)}+2g(D)
 -2b-n'-2q_G+2g(G)\nonumber\\
&\hspace{20mm}-2a-q_1-q_2+g(F_1)+g(F_2).
\label{eq:B6}
\end{align}

\begin{proposition}\label{prop:balance}
With the notation above, $\mathcal B_6=0$.
This relation follows from topological Lefschetz, the point decomposition and
Riemann--Hurwitz for the two residual actions.
\end{proposition}

\begin{proof}
Write
\begin{align*}
\mathcal R_G={}&4q_G-2g(G)-2g(D)+p_{(3,4)}-2k+4b+2\ell,\\
\mathcal R_F={}&6(q_1+q_2)-2g(F_1)-2g(F_2)-4g(D)
 +2p_{(2,5)}+2p_{(3,4)}-4N+12a+4\ell.
\end{align*}
By \cref{eq:RHG,eq:RHF}, both forms vanish.  Expanding the following
expression, we obtain
\[
\mathcal B_6=\frac12\mathcal P_6+\frac16\mathcal T_1
 +\frac12\mathcal T_2+\frac13\mathcal T_3
 -\frac12\mathcal R_G-\frac16\mathcal R_F.
\]
All terms on the right are zero.  This proves the proposition.
\end{proof}

Consider also the linear form
\begin{align}
\mathcal C_6:={}&-n'-3+\frac32r-6\ell-2p_{(2,5)}-3p_{(3,4)}+6g(D)
 +2k-6b-6q_G+4g(G)\nonumber\\
&+\frac32N-6a-3q_1-3q_2+\frac32g(F_1)+\frac32g(F_2).
\label{eq:C6}
\end{align}

\begin{proposition}\label{prop:dependence}
The following identity of linear forms holds:
\begin{equation}\label{eq:dependency}
\mathcal C_6=\frac14\mathcal D_6-\mathcal P_6-\mathcal T_2
 -\frac34\mathcal T_3+3\mathcal B_6.
\end{equation}
Thus $\mathcal C_6=0$ follows from the dimension formula, the point
decomposition, topological Lefschetz for $\sigma^2,\sigma^3$, and the balance
relation.
\end{proposition}

\begin{proof}
Expanding the right-hand side of \cref{eq:dependency} and collecting
coefficients gives \cref{eq:C6}.
\end{proof}

\medskip
The linear form on the left-hand side of
\[
-2\alpha+10+N-r-g(F_1)-g(F_2)=0
\]
equals $-\tfrac12(\mathcal D_6+\mathcal T_3)$ and hence gives no independent
relation.

\subsection{Explicit \texorpdfstring{order-$6$}{order-6} models}

\begin{example}\label{ex:generic6}
Let $f_{12}$ have twelve simple zeros and consider
\[
S_f:\quad y^2=x^3+f_{12}(t),\qquad
\sigma(x,y,t)=(\zeta_3x,-y,t).
\]
The action on the holomorphic two-form is by a primitive sixth root.  The fixed
curves are the zero section, a bisection fixed by $\sigma^2$, and a trisection
fixed by $\sigma^3$.  Their ramification over the twelve zeros of $f_{12}$ gives
\begin{gather*}
(r,m,\alpha,\beta)=(2,10,0,0),\qquad
(n,n',k,g(G))=(0,0,2,5),\\
(p_{(3,4)},p_{(2,5)},\ell,N,b,a,g(F_1)+g(F_2))
=(12,0,1,2,0,0,10).
\end{gather*}
All quotient genera in \cref{eq:B6} are zero, and substitution gives
$\mathcal B_6=0$.  The associated Calabi--Yau threefold has Hodge numbers
$(29,29)$~\cite[Table~5]{CattaneoGarbagnati}.
\end{example}

\medskip
\begin{example}\label{ex:special6}
For $\lambda\notin\{0,\pm1\}$ take
\[
f_{12}(t)=t^4(t-1)^3(t+1)^3(t-\lambda)^2.
\]
Resolving the four reducible fibers gives
\begin{gather*}
(r,m,\alpha,\beta)=(11,2,2,3),\qquad
(n,n',k,g(G))=(6,1,3,0),\\
(p_{(3,4)},p_{(2,5)},\ell,N,b,a,g(F_1)+g(F_2))
=(4,4,1,5,0,0,0).
\end{gather*}
\Cref{eq:D6,eq:P6,eq:T1,eq:T2,eq:T3,eq:H6,eq:B6} all vanish.
The associated threefold has Hodge numbers $(55,1)$
\cite[Example~7.9 and Table~5]{CattaneoGarbagnati}.
\end{example}

The preceding two examples use the same construction.  It also gives all the
rational-curve cases in the final table.  We shall write
$(a_1,a_2);kR$ if the fixed locus consists of $k$ rational curves,
$a_1$ points of type $(2,5)$ and $a_2$ points of type $(3,4)$.

\begin{proposition}\label{prop:weierstrass-six}
Let $P(s,t)$ be a binary form of degree $12$ whose roots have multiplicities
$1$, $2$, $4$ or $5$, and let $X_P$ be the minimal model of
\begin{equation}\label{eq:weierstrass-six}
 y^2=x^3+P(s,t).
\end{equation}
Then $X_P$ is a $K3$ surface and
\[
 \sigma_0(x,y;[s:t])=(\zeta_3x,-y;[s:t])
\]
is purely non-symplectic of order $6$.  If $s_i$ denotes the number of roots of
multiplicity $i$, then
\begin{equation}\label{eq:weierstrass-fixed-data}
 p_{(3,4)}=s_1+2s_4+3s_5,\qquad
 p_{(2,5)}=s_2+s_4+4s_5,\qquad
 \ell=1+s_5.
\end{equation}
In particular, the following fixed loci are obtained.  In the second column
$1^{s_1}2^{s_2}4^{s_4}5^{s_5}$ denotes the multiplicities of the roots of
$P$.
\begin{center}
\begin{tabular}{@{}ll@{}}
\toprule
$\Fix(\sigma_0)$ & multiplicities of the roots of $P$\\
\midrule
$(0,12);R$ & $1^{12}$\\
$(1,10);R$ & $1^{10}2$\\
$(2,8);R$ & $1^8 2^2$\\
$(3,6);R$ & $1^6 2^3$\\
$(4,4);R$ & $1^4 2^4$\\
$(5,2);R$ & $1^2 2^5$\\
$(4,10);2R$ & $1^7 5$\\
$(5,8);2R$ & $1^5 2 5$\\
$(6,6);2R$ & $1^3 2^2 5$\\
$(7,4);2R$ & $1 2^3 5$\\
$(8,8);3R$ & $1^2 5^2$\\
$(9,6);3R$ & $2 5^2$\\
$(6,0);R$ & $2^6$\\
\bottomrule
\end{tabular}
\end{center}
For the multiplicities $\mu_1,\ldots,\mu_q$ in the table we may take
\begin{equation}\label{eq:explicit-root-form}
 P(s,t)=\prod_{j=1}^q(t-a_j s)^{\mu_j},
 \qquad a_i\ne a_j,
\end{equation}
where $\sum\mu_j=12$.
\end{proposition}

\begin{proof}
The fibers over roots of multiplicity $1$, $2$, $4$ and $5$ have Kodaira types
$II$, $IV$, $IV^*$ and $II^*$, respectively.  The zero section is fixed by
$\sigma_0$.  On the four fibers the remaining contributions are, in the same
order,
\[
 P_{(3,4)},\qquad P_{(2,5)},\qquad
 2P_{(3,4)}+P_{(2,5)},\qquad
 3P_{(3,4)}+4P_{(2,5)}+R.
\]
This is the local calculation of Dillies
\cite[Lemma~7.4]{Dillies6}; it gives
\cref{eq:weierstrass-fixed-data}.  The sum of the Euler numbers of the fibers is
$24$.  The minimal Weierstrass fibration has fundamental line bundle
$\mathcal O_{\PP^1}(2)$, so its canonical bundle is trivial and
$H^1(X_P,\mathcal O_{X_P})=0$.  Thus $X_P$ is a $K3$ surface.  Finally,
$dt\wedge dx/y$ is multiplied by $-\zeta_3$, a primitive sixth root of unity.
\end{proof}

For the last line of the table, let $q_6$ be a binary form with six distinct
roots and put $P=q_6^2$.  The surface has six
fibers of type $IV$.  The automorphism $\sigma_0$ fixes the zero section and one
point of type $(2,5)$ on every reducible fiber.  Moreover
\begin{equation}\label{eq:square-model-fixes}
 \Fix(\sigma_0)=R\sqcup6P_{(2,5)},\quad
 \Fix(\sigma_0^2)=3R\sqcup6P,\quad
 \Fix(\sigma_0^3)=R\sqcup C_4.
\end{equation}
Here $C_4$ is the normalization of $u^3=q_6(s,t)$.  This is the class
\texttt{0.6.3.15} in the list of Brandhorst and Hofmann.

A second lift of the order-$3$ action gives examples without fixed curves.

\begin{proposition}\label{prop:base-involution-six}
Consider the three binary forms
\[
 P_{00}=s^{12}+t^{12},\qquad
 P_{20}=t^2(s^{10}+t^{10}),\qquad
 P_{22}=s^2t^2(s^8+t^8).
\]
On the minimal model of $y^2=x^3+P_{ij}(s,t)$ let
\begin{equation}\label{eq:base-involution-action}
 \sigma_1(x,y;[s:t])=(\zeta_3x,y;[s:-t]).
\end{equation}
Then $\sigma_1$ is purely non-symplectic of order $6$ and its fixed locus is,
respectively,
\[
 (0,6);-,\qquad (1,4);-,\qquad (2,2);-.
\]
\end{proposition}

\begin{proof}
The fixed locus is contained in the fibers over $[1:0]$ and $[0:1]$.  For
$P_{00}$ both fibers are smooth.  For $P_{20}$ their types are $IV$ and $I_0$,
and for $P_{22}$ they are both of type $IV$.  A preserved smooth fiber contains
three points of type $(3,4)$, while a preserved fiber of type $IV$ contains one
point of each type~\cite[Lemma~7.5]{Dillies6}.  This gives the three fixed loci.
The form $dt\wedge dx/y$ is multiplied by $-\zeta_3$.
\end{proof}

\medskip
\begin{example}
\label{ex:complete-intersection-six}
Let $X\subset\PP^4$ be given by
\begin{align*}
 x_0^2-x_1^2-x_2x_3+x_2x_4&=0,\\
 x_0^3+x_2^3+x_3^3+x_0x_2x_4+x_3x_4^2&=0.
\end{align*}
The Jacobian criterion shows that $X$ is a smooth $K3$ surface.  Consider the
projective transformation
\begin{equation}\label{eq:complete-intersection-action}
 \sigma(x_0,x_1,x_2,x_3,x_4)
 =(x_0,-x_1,\zeta_3x_2,\zeta_3^2x_3,\zeta_3^2x_4)
\end{equation}
\par\noindent
The transformation preserves both equations.  Its determinant is $\zeta_6$,
so it multiplies the residue form by a primitive sixth root.  The fixed locus consists of the three
points on the line $x_0=x_1=x_2=0$ cut out by
$x_3(x_3^2+x_4^2)=0$.  Since no curve is fixed, the holomorphic Lefschetz
formula shows that they are all of type $(2,5)$.  Furthermore
\[
 \Fix(\sigma^2)=3P,\qquad \Fix(\sigma^3)=C_4,
\]
where $C_4$ is the smooth complete intersection in the hyperplane $x_1=0$.
Thus we obtain the class \texttt{0.6.4.26} and the generator locus
$(3,0);-$.
\end{example}

\medskip
\begin{example}\label{ex:genus2-six}
Let $F_6(x_0,x_1)$ be square-free and $c\ne0$.  The double plane
\begin{equation}\label{eq:genus2-model}
 X_{F,c}:\quad w^2=F_6(x_0,x_1)+c x_2^6
 \subset\PP(1,1,1,3)
\end{equation}
is a smooth $K3$ surface for general coefficients.  It carries
\[
 \sigma(x_0,x_1,x_2;w)=(x_0,x_1,\zeta_6x_2;w).
\]
The residue two-form is multiplied by $\zeta_6$.  The line $x_2=0$ lifts to
the genus-$2$ curve $C_2:w^2=F_6(x_0,x_1)$ and is fixed pointwise.  The two
points over $[0:0:1]$ are exchanged by $\sigma$, fixed separately by
$\sigma^2$, and exchanged again by $\sigma^3$.  Consequently
\begin{equation}\label{eq:genus2-fixes}
 \Fix(\sigma)=C_2,\qquad
 \Fix(\sigma^2)=C_2\sqcup2P,\qquad
 \Fix(\sigma^3)=C_2.
\end{equation}
Thus $(r,m,\alpha,\beta)=(1,4,4,5)$.  This is class \texttt{0.6.3.37} in the
tables of~\cite{BrandhorstHofmannData}.
\end{example}

\medskip
\begin{example}\label{ex:elliptic-six}
The minimal model of
\begin{equation}\label{eq:dillies-elliptic-model}
 y^2=x^3+(t^6-1)^2,
 \qquad \rho(x,y,t)=(x,y,\zeta_6t),
\end{equation}
is a $K3$ surface with a purely non-symplectic automorphism of order $6$.  The
fiber $E_0$ over $t=0$ is fixed pointwise.  Dillies' local calculation gives
\cite[Proposition~6.2]{Dillies6}
\[
 \Fix(\rho)=\Fix(\rho^2)=E_0\sqcup3P_{(2,5)},\qquad
 \Fix(\rho^3)=E_0\sqcup E_\infty.
\]
It realizes the class \texttt{0.6.3.36} and the generator locus $(3,0);C_1$.
\end{example}

\medskip
\begin{example}\label{ex:missing-six}
Let $a_0,a_6,b_0,b_6\in\CC$ belong to the non-empty Zariski-open set on which
$a_6\ne0$, $4a_0^3+27b_0^2\ne0$, and the finite discriminant has eighteen
simple zeros.  This set contains $(a_0,a_6,b_0,b_6)=(0,1,1,0)$, for which
the discriminant is $-16(4t^{18}+27)$.  Let $X$ be the minimal smooth model of
\begin{equation}\label{eq:missing-model}
 y^2=x^3+(a_0+a_6t^6)x+(b_0+b_6t^6),\qquad
 \sigma(x,y,t)=(x,y,\zeta_6t).
\end{equation}
In homogeneous base coordinates the two Weierstrass coefficients are
\[
 a_0s^8+a_6s^2t^6\in H^0(\PP^1,\mathcal O(8)),\qquad
 b_0s^{12}+b_6s^6t^6\in H^0(\PP^1,\mathcal O(12)).
\]
Thus the relatively minimal resolution has fundamental line bundle
$\mathcal O_{\PP^1}(2)$ and trivial canonical bundle; the genericity assumptions give
eighteen $I_1$ fibers and the $I_0^*$ fiber described below, so its Euler number
is $18+6=24$.  Since the base is $\PP^1$ and the fundamental line bundle is
$\mathcal O_{\PP^1}(2)$, we also have $H^1(X,\mathcal O_X)=0$.
It is therefore a $K3$ surface.  Moreover
$\omega_X=dt\wedge dx/y$ on the affine chart, and hence
$\sigma^*\omega_X=\zeta_6\omega_X$.

The fiber $E_0$ over $t=0$ is smooth and fixed pointwise.  A point fixed by
$\sigma$ must lie over one of the two fixed points of the base action.  At
infinity put
$u=t^{-1}$, $X=u^4x$, $Y=u^6y$.  The equation and action become
\begin{align*}
Y^2&=X^3+u^2(a_6+a_0u^6)X+u^6(b_6+b_0u^6),\\
(u,X,Y)&\longmapsto(\zeta_6^{-1}u,\zeta_6^2X,Y).
\end{align*}
Here $\operatorname{ord}_0(A)=2$, $\operatorname{ord}_0(B)\geq6$, and
$\operatorname{ord}_0(\Delta)=6$; Tate's algorithm
therefore gives a fiber of type $I_0^*$.  We now resolve this fiber and compute the action on the exceptional curves
to find the four isolated fixed points.  In the $u$-chart of the
first blow-up, set $X=us$ and $Y=uv$.  The strict transform and the lifted action
are
\begin{align}
 v^2&=u(s^3+a_6s)+b_6u^4+a_0u^7s+b_0u^{10},\label{eq:missing-blowup}\\
 (u,s,v)&\longmapsto(\zeta_6^{-1}u,-s,\zeta_6v).\label{eq:missing-action}
\end{align}
The exceptional line is $u=v=0$.  It contains the three ordinary double
points
\[
 s(s^2+a_6)=0.
\]
The two non-zero points are exchanged, while the point $s=0$ is preserved.
Resolving the latter node, locally $v^2=us$, gives in its two standard charts
tangent eigenvalues $(\zeta_6^5,\zeta_6^2)$ and
$(\zeta_6^3,\zeta_6^4)$.  Thus this arm of the $\widetilde D_4$ fiber contributes
one point of type $(2,5)$ and one of type $(3,4)$.

We find the fourth arm in the $X$-chart of the first blow-up.  Writing
$u=Xq$ and $Y=Xw$ gives
\[
 w^2=X(1+a_6q^2)+b_6X^4q^6+a_0X^7q^8+b_0X^{10}q^{12},
 \qquad(q,w)\longmapsto(-q,\zeta_6^4w).
\]
At its intersection with the central component the tangent type is $(3,4)$.
The other fixed point of the induced order-$3$ action on that rational arm has
tangent weight $2$ along the arm; since the product of the two tangent
eigenvalues is $\zeta_6$, its type is $(2,5)$.  The remaining two arms are
exchanged.  Thus the local calculation gives
\begin{equation}\label{eq:missing-fixes}
\begin{array}{c|c}
g&\Fix(g)\\ \hline
\sigma&E_0\sqcup2P_{(2,5)}\sqcup2P_{(3,4)}\\
\sigma^2&E_0\sqcup R\sqcup4P\\
\sigma^3&E_0\sqcup2R.
\end{array}
\end{equation}
Indeed, $\sigma^2$ fixes the central component and acts non-trivially on each
of the four arms, producing one further fixed point on every arm.  The two arms
preserved by $\sigma$ are fixed pointwise by $\sigma^3$, whereas the other two
are exchanged.  The holomorphic Lefschetz equation
$p_{(3,4)}+2p_{(2,5)}=6$ gives the split $(2,2)$.  The Euler characteristics in
\cref{eq:missing-fixes} give
$(r,m,\alpha,\beta)=(6,3,3,4)$ by
\cref{eq:coh-reconstruction6}.  This realizes the class
\texttt{0.6.2.29}, whose generator locus $E_0\sqcup4P$ was absent from
Dillies' list~\cite{Dillies6,BrandhorstHofmann}.  Here $n'=1$, $q_G=1$, and
$q_1+q_2=1$, and $\mathcal B_6=0$ follows as well.
\end{example}

\medskip
\begin{corollary}
\label{cor:twenty-explicit}
Every one of the $20$ possible fixed loci of a purely non-symplectic
automorphism of order $6$ in \cref{thm:numerical-classification} is realized by
an explicit equation.
\end{corollary}

\begin{proof}
The thirteen loci containing only rational curves and isolated points are given
by Proposition~\ref{prop:weierstrass-six}.  The three loci $(0,6);-$, $(1,4);-$
and $(2,2);-$ are given by Proposition~\ref{prop:base-involution-six}, while
Example~\ref{ex:complete-intersection-six} gives $(3,0);-$.  Finally,
Examples~\ref{ex:genus2-six}, \ref{ex:elliptic-six} and
\ref{ex:missing-six} give
$(0,0);C_2$, $(3,0);C_1$ and $(2,2);C_1$, respectively.
\end{proof}

\section{Enumeration in order six}
\label{sec:sarti-dillies}

Put
\[
 \tau=\sigma^2,\qquad \iota=\sigma^3.
\]
Thus $\tau$ is a non-symplectic automorphism of order $3$, $\iota$ is a
non-symplectic involution, and the two commute.  Following
Sarti and Dillies, we study the residual actions on these two fixed loci.  We
use the $24$ order-$3$ types of
Artebani--Sarti~\cite{ArtebaniSarti3} and the $75$ deformation types of
non-symplectic involutions (with $66$ fixed-locus topologies) from Nikulin's
classification~\cite{NikulinInvolutions}.

\medskip
\begin{lemma}\label{lem:residual-actions}
Let $C$ be a component of $\Fix(\tau)$.  Under the residual involution induced
by $\sigma$, precisely one of the following occurs:
\begin{enumerate}[label=(\alph*)]
\item $C$ is fixed pointwise and hence belongs to $\Fix(\sigma)$;
\item $C\simeq\PP^1$ is invariant non-trivially and contributes two points of
type $(3,4)$;
\item $C$ has genus $h>0$, is invariant but not fixed pointwise, with quotient
genus $q$, and contributes
$2h+2-4q$ points of type $(3,4)$;
\item $C$ belongs to an exchanged pair.
\end{enumerate}
Similarly, a component $F\subset\Fix(\iota)$ is fixed pointwise, belongs to a
triple, or carries a non-trivial residual order-$3$ action.  In the last case it
contributes two fixed points if $F$ is rational and
\begin{equation}\label{eq:RH-local-three}
 h+2-3q
\end{equation}
fixed points if $g(F)=h>0$ and $g(F/\sigma)=q$.  Finally, if $n$ is the number
of isolated points of $\tau$, then $n=p_{(2,5)}+2n'$.
\end{lemma}

\begin{proof}
Only orbit lengths $1,2$ occur on $\Fix(\tau)$ and only $1,3$ on
$\Fix(\iota)$.  A non-trivial automorphism of a rational curve has two fixed
points.  For a positive-genus curve the two formulas are Riemann--Hurwitz:
\[
 2h-2=2(2q-2)+f,\qquad
 2h-2=3(2q-2)+2f.
\]
Squaring the local action at a point of type $(2,5)$ gives an isolated fixed point of
$\tau$, while a point of type $(3,4)$ lies on a $\tau$-fixed curve.  This proves the last
statement and identifies the tangent types.
\end{proof}

\medskip
Using the lemma, we enumerate the possible fixed loci.  For each pair of
lower-order types we choose the pointwise-fixed components, the paired or
tripled components and the quotient genera.  We then impose Riemann--Hurwitz,
\cref{eq:D6,eq:P6,eq:T1,eq:T2,eq:T3,eq:H6,eq:B6}, non-negativity, containment
of fixed loci and the Hodge-index intersection bound.  The cohomological
multiplicities are reconstructed from the three Euler characteristics by
\cref{eq:coh-reconstruction6}.

\subsection{The local lattice sieve}

Let $L=H^2(S,\ZZ)$, $L_2=L^{\sigma^2}$ and $L_3=L^{\sigma^3}$.  For an
even lattice $M$ we denote its discriminant form by $q_M$ and its $p$-primary
part by $q_M(p)$.  Put
\[
 C=\ker(\sigma^*-1),\quad A=\ker(\sigma^*+1),\quad
 B=\ker\Phi_3(\sigma^*),\quad D=\ker\Phi_6(\sigma^*).
\]
Their signatures are
\begin{equation}\label{eq:cyclotomic-signatures}
 (1,r-1),\qquad (0,\beta),\qquad (0,2\alpha),\qquad
 (2,2m-2),
\end{equation}
respectively.  The genera of $L_2$ and $L_3$ are determined by the
order-$3$ and order-$2$ classifications.

\begin{lemma}\label{lem:local-splitting}
There are non-degenerate finite quadratic forms $H_2,K_2$ on elementary
$2$-groups and $H_3,K_3$ on elementary $3$-groups such that
\begin{equation}\label{eq:lower-power-splitting}
 q_{L_3}=H_2\oplus K_2,\qquad q_{L_2}=H_3\oplus K_3,
\end{equation}
and
\begin{equation}\label{eq:cyclotomic-discriminant-forms}
\begin{alignedat}{2}
q_C&=H_2\oplus H_3,\qquad& q_A&=(-H_2)\oplus K_3,\\
q_B&=K_2\oplus(-H_3),& q_D&=(-K_2)\oplus(-K_3).
\end{alignedat}
\end{equation}
\end{lemma}

\begin{proof}
The only non-trivial resultants between
$\Phi_1,\Phi_2,\Phi_3,\Phi_6$ are
\[
 |\Res(\Phi_1,\Phi_2)|=2,\quad
 |\Res(\Phi_1,\Phi_3)|=|\Res(\Phi_2,\Phi_6)|=3,
 \quad |\Res(\Phi_3,\Phi_6)|=4.
\]
Consequently all gluings are supported at $2$ and $3$.  The sum
$C\oplus B$ has no $2$-primary gluing, while $C\oplus A$ has no
$3$-primary gluing.  This gives \cref{eq:lower-power-splitting}.
Unimodularity of $L$ pairs the remaining primary forms with the opposite
sign and gives \cref{eq:cyclotomic-discriminant-forms}; see
\cite[Sections~2.2 and~4.2]{BrandhorstHofmann} and
\cite[Section~1]{NikulinForms}.
\end{proof}

Let $v_2$ be the quadratic form on $(\ZZ/2\ZZ)^2$ with matrix
\[
 \begin{pmatrix}1&1/2\\1/2&1\end{pmatrix}
 \quad\text{in }\QQ/2\ZZ.
\]

\begin{lemma}\label{lem:eisenstein-obstructions}
Every splitting in Lemma~\ref{lem:local-splitting} which comes from an
order-$6$ isometry satisfies
\begin{equation}\label{eq:local-2-condition}
 K_2\simeq v_2^{\oplus t}
\end{equation}
for some $t$, and
\begin{equation}\label{eq:local-3-condition}
 \ell(H_3)\leq\alpha,\quad \ell(H_3)\equiv\alpha\pmod2,
 \qquad
 \ell(K_3)\leq m,\quad \ell(K_3)\equiv m\pmod2.
\end{equation}
Moreover $B_\QQ$ and $D_\QQ$ are trace forms of Hermitian spaces over
$\QQ(\zeta_3)$.
\end{lemma}

\begin{proof}
The lattices $B$ and $D$ are modules over
$\ZZ[\zeta_3]$.  Since $2$ is inert, reduction modulo $2$ makes their
$2$-elementary discriminant groups vector spaces over $\FF_4$.  The unique
non-degenerate Hermitian form on a one-dimensional $\FF_4$-space has
underlying quadratic form $v_2$.  Orthogonal diagonalization gives
\cref{eq:local-2-condition}.

The prime $3$ is ramified in $\ZZ[\zeta_3]$.  The equivariant gluing to $C$ makes the action on $H_3$ trivial, while
the gluing to $A$ makes the action on $K_3$ equal to $-1$.
Thus the corresponding Eisenstein discriminant quotients are killed by
$1-\zeta_3$; compare \cite[Proposition~4.7]{BrandhorstHofmann}.
Such a quotient of an Eisenstein lattice of rank $e$ is generated by at most
$e$ elements.  The determinant of its trace lattice has square class
$3^e$.  Applying this to the Eisenstein ranks $\alpha$ and $m$ gives
\cref{eq:local-3-condition}.  The last assertion follows from the correspondence between
cyclotomic isometries and Hermitian forms
\cite[Section~4.1]{BrandhorstHofmann}.
\end{proof}

We check the rational Hermitian condition without constructing an integral
Hermitian lattice.  Set
\[
 A_2=\begin{pmatrix}2&-1\\-1&2\end{pmatrix}.
\]
After diagonalizing a Hermitian space, its trace form is an orthogonal sum of
blocks $aA_2$.  The trace lattice is unimodular away from $2$ and $3$.  Therefore the
Hermitian determinant has even valuation at every inert prime other than $2$.
Since $3=N_{\QQ(\zeta_3)/\QQ}(1-\zeta_3)$, its norm class has a
representative in $\{1,2,-1,-2\}$.
For a sum of trace blocks the local Hasse invariant depends only on the rank
and the product of the coefficients modulo norms: the variable factor is
$(\prod_i a_i,-3)_p$.  Consequently the Hasse invariants at $2$ and $3$ of
$B_\QQ$ must agree with those of a sum of $\alpha$ blocks with coefficients
$-1$ or $-2$.  For $D_\QQ$ one takes one coefficient in $\{1,2\}$ and $m-1$
coefficients in $\{-1,-2\}$.  Together with the determinant condition in
\cref{eq:local-3-condition}, the Hasse--Minkowski theorem makes this a finite
test.  We use this condition below; see
\cite[Section~4.1]{BrandhorstHofmann}.

For every splitting we also require that the four forms in
\cref{eq:cyclotomic-discriminant-forms}, with the signatures in
\cref{eq:cyclotomic-signatures}, are discriminant forms of even integral
lattices.  This is decided by the local existence criterion of
\cite[Theorem~1.10.1]{NikulinForms}.  Finally, over $\QQ$ we require
\begin{equation}\label{eq:rational-compatibility}
 C_\QQ\perp A_\QQ\simeq (L_2)_\QQ,\qquad
 C_\QQ\perp B_\QQ\simeq (L_3)_\QQ.
\end{equation}
We check \cref{eq:rational-compatibility} by comparing dimensions,
determinants and local Hasse invariants.

\begin{lemma}\label{lem:root-free-bound}
Let $M$ be either $-A$ or $-B$ and let $s=\rk M>0$.  Then $M$ is positive
definite, has minimum at least $4$, and
\begin{equation}\label{eq:root-free-determinant-bound}
 \det M\geq d_s,
\end{equation}
where the values needed here are
\[
\begin{array}{c|rrrrrrrrrrrr}
s&1&2&3&4&5&6&7&8&9&10&12&14\\ \hline
d_s&4&12&32&61&106&159&217&256&298&317&281&177.
\end{array}
\]
\end{lemma}

\begin{proof}
An invariant ample class belongs to $C$.  A vector of square $-2$ in $A$ or
$B$ would be orthogonal to this class, which is impossible by
Riemann--Roch.  Hence $M$ is root-free.  The lattice packing obtained from
$M$ has centre density at least $1/\sqrt{\det M}$.  If $\Delta_s$ is an upper
bound for the packing density in $\RR^s$ and
$V_s=\pi^{s/2}/\Gamma(s/2+1)$, then
\[
 \frac1{\sqrt{\det M}}\leq\frac{\Delta_s}{V_s},
 \qquad
 \det M\geq\left(\frac{V_s}{\Delta_s}\right)^2.
\]
The displayed integers are the ceilings obtained from the rigorous
three-point bounds of Cohn, de Laat and Salmon
\cite[Table~1.1]{CohnDeLaatSalmon}, whose displayed upper bounds are rounded
upwards.  In
dimensions $1,2,3$ and $8$ the exact optimal densities are used.  The
classical linear-programming bounds of~\cite{CohnElkies} would leave three
additional candidates.
\end{proof}

\begin{proposition}
\label{thm:geometric-classification}
The residual-action enumeration gives $817$ triples.  The necessary conditions
on the four cyclotomic sublattices leave $204$ of them.  All the $142$ fixed-locus
triples in the Brandhorst--Hofmann list occur among these $204$ triples.
\end{proposition}

\begin{proof}
We use the possibilities from Lemma~\ref{lem:residual-actions} and the
necessary equations above.  For each of the $24$ order-$3$ types and $66$ involution fixed-locus
topologies, we choose the pointwise-fixed components and then the components
in orbits of lengths $2$ and $3$.  The two Riemann--Hurwitz formulas determine
the possible quotient genera.  Removing repetitions gives $817$ triples.

For every triple we enumerate the non-degenerate splittings in
\cref{eq:lower-power-splitting}.  We apply the integral genus test,
\cref{eq:local-2-condition,eq:local-3-condition}, the rational compatibility
\cref{eq:rational-compatibility}, the Hermitian Hasse test and
Lemma~\ref{lem:root-free-bound}.  All sets involved are finite: an elementary
$3$-form is determined by its length and determinant, and the elementary
$2$-forms are orthogonal sums of the standard one- and two-dimensional forms.
Exactly $204$ triples admit at least one splitting.
\end{proof}

The following table gives the counts for each fixed locus of the generator.
The first column uses the preceding notation.
\begin{center}
\scriptsize
\begin{longtable}{@{}lrrr@{}}
\toprule
$\Fix(\sigma)$ & geometric & local lattice & realized\\
\midrule
\endfirsthead
\toprule
$\Fix(\sigma)$ & geometric & local lattice & realized\\
\midrule
\endhead
$(0,0);C_2$ & 7&1&1\\
$(0,0);R+C_3$ & 4&0&0\\
$(0,12);R$ & 11&1&1\\
$(0,6);-$ & 149&32&17\\
$(0,6);C_1$ & 3&0&0\\
$(1,10);R$ & 23&3&3\\
$(1,4);-$ & 125&32&25\\
$(2,2);-$ & 126&36&26\\
$(2,2);C_1$ & 6&1&1\\
$(2,2);R+C_2$ & 3&0&0\\
$(2,2);2R+C_3$ & 1&0&0\\
$(2,8);R$ & 43&7&6\\
$(2,8);R+C_1$ & 1&0&0\\
$(3,0);-$ & 54&26&15\\
$(3,0);C_1$ & 3&1&1\\
$(3,6);R$ & 55&13&12\\
$(4,10);2R$ & 4&1&1\\
$(4,4);R$ & 69&15&14\\
$(4,4);R+C_1$ & 3&0&0\\
$(5,2);R$ & 50&13&6\\
$(5,8);2R$ & 7&3&3\\
$(6,0);R$ & 28&9&1\\
$(6,0);R+C_1$ & 3&0&0\\
$(6,6);2R$ & 13&5&5\\
$(6,6);2R+C_1$ & 1&0&0\\
$(7,4);2R$ & 10&2&2\\
$(8,2);2R$ & 10&1&0\\
$(8,2);2R+C_1$ & 1&0&0\\
$(8,8);3R$ & 1&1&1\\
$(9,0);2R$ & 2&0&0\\
$(9,6);3R$ & 1&1&1\\
\midrule
total & 817&204&142\\
\bottomrule
\end{longtable}
\end{center}

The realized column is obtained by grouping the $150$ order-$6$ rows of
Brandhorst and Hofmann by the three fixed loci.  There are $134$ singleton
groups and eight groups of size two, listed in Section~\ref{sec:class-program};
hence there are $134+8=142$ groups.  Every one of them passes the local tests.

\medskip
For the $204$ remaining triples, a classification still requires compatible
actions on representatives of the definite genera and an equivariant primitive
extension in the $K3$ lattice.  We therefore use the complete classification
of Brandhorst and Hofmann in the next section.

\subsection{Relation with the lower-order models}

The square of the automorphism in Proposition~\ref{prop:weierstrass-six} is the
Artebani--Sarti action
\[
 (x,y,t)\longmapsto(\zeta_3^2x,y,t)
\]
on $y^2=x^3+p_{12}(t)$.  For fixed loci with more than one curve, the
multiplicities of the roots of $p_{12}$ are described in
\cite[Proposition~4.2]{ArtebaniSarti3}.  Thus
Proposition~\ref{prop:weierstrass-six} gives the order-$6$ lift and its action
on every exceptional component.  The cases with at most one curve are
represented by the $(2,3)$ complete intersections, quartics and double sextics
of Artebani and Sarti~\cite[Propositions~4.7, 4.9 and 4.11]{ArtebaniSarti3};
Examples~\ref{ex:complete-intersection-six} and \ref{ex:genus2-six} give two
such lifts.

For the cube $\iota$, one may use the usual double-plane models for
non-symplectic involutions.  An irreducible nodal sextic whose normalization has
genus $g$ gives the component $C_g$, while a transverse union of a line and a
smooth quintic gives $R\sqcup C_6$.  The models for the two powers must induce
compatible actions on the same $K3$ lattice.  This global condition is not checked in
Proposition~\ref{thm:geometric-classification}.  If
commuting representatives $\tau$ and $\iota$ have been found, then
\begin{equation}\label{eq:lift-generator}
 \sigma=\tau^2\circ\iota,
 \qquad \sigma^2=\tau,
 \qquad \sigma^3=\iota.
\end{equation}
Different projective models for the two powers may be related by a change of
polarization, as in the constructions of Artebani--Sarti and
Dillies~\cite[Sections~6--7]{Dillies6}.

\section{Applications}

\subsection{Numerical fixed-locus types in order six}\label{sec:class-program}

Brandhorst and Hofmann find exactly $150$ deformation classes of purely
non-symplectic order-$6$ actions and list the fixed loci of $\sigma$,
$\sigma^2$, and $\sigma^3$
\cite{BrandhorstHofmann,BrandhorstHofmannData}.  Their result also contains the
case missed in Dillies' list, namely a genus-$1$ curve and four isolated fixed
points for $\sigma$ (class \texttt{0.6.2.29} in their tables).

We apply the preceding formulas to their list.  We write a fixed locus as
\[
 ((a_1,a_2),k,l,g),
\]
where $a_1=p_{(2,5)}$, $a_2=p_{(3,4)}$, $k$ is the number of rational curves,
and $l$ is the number of positive-genus curves, each of genus $g$; for the
proper powers the point tuple has length one or zero.  The symbols $k,l$ here have different meanings from the curve counts
$k,\ell$ of Section~5.  Set
$\varepsilon_j=\eul(\Fix(\sigma^j))$.  The dimension
and topological Lefschetz relations have the unique solution
\begin{align}\label{eq:coh-reconstruction6}
r&=\frac{2\varepsilon_1+2\varepsilon_2+\varepsilon_3+12}{6},&
m&=\frac{\varepsilon_1-\varepsilon_2-\varepsilon_3+24}{6},\nonumber\\
\alpha&=\frac{-\varepsilon_1-\varepsilon_2+\varepsilon_3+24}{6},&
\beta&=\frac{-2\varepsilon_1+2\varepsilon_2-\varepsilon_3+24}{6}.
\end{align}
\par\noindent
The point-orbit relation then recovers
\begin{equation}\label{eq:nprime-reconstruction}
n'=\frac{n-p_{(2,5)}}2.
\end{equation}
Thus the numerical types of the three fixed loci determine the cohomological
multiplicities.  For a candidate, the right-hand sides must be non-negative integers.

In the next theorem, a fixed locus means its point counts and the genera of
its connected components.

\begin{theorem}
\label{thm:numerical-classification}
The $150$ purely non-symplectic order-$6$ deformation classes of Brandhorst and
Hofmann give exactly $142$ sets of numerical data
\[
 \mathfrak N(\sigma)=
 \bigl(\Fix(\sigma),\Fix(\sigma^2),\Fix(\sigma^3);
       r,m,\alpha,\beta,n'\bigr).
\]
Among them, $134$ occur for one deformation class and eight occur for
two classes.  The generator has $20$ possible fixed loci.  There are $102$
different sequences $\{\eul(Y_{6,n})\}_{n\geq1}$.
\end{theorem}

\begin{proof}
Apply \cref{eq:coh-reconstruction6,eq:nprime-reconstruction} to the
Brandhorst--Hofmann list.  This gives $150$ rows and $142$ distinct sets of data
with multiplicity distribution
$134\cdot1+8\cdot2$, and $20$ generator loci.  There are $102$ distinct triples
$(\varepsilon_1,\varepsilon_2,\varepsilon_3)$.  By
Proposition~\ref{prop:euler-determination}, these give exactly $102$ distinct
Euler sequences.  The script \texttt{analyze\_bh\_order6.py} groups the rows of version~2.1 of
the database.
\end{proof}

The following table gives the possible fixed loci of the generator.  We write
$(a_1,a_2);kR+lC_g$ for the preceding tuple.
\begin{center}
\small
\begin{longtable}{@{}lrrr@{}}
\toprule
$\Fix(\sigma)$ & deformation classes & numerical types & stringy types\\
\midrule
\endfirsthead
\toprule
$\Fix(\sigma)$ & deformation classes & numerical types & stringy types\\
\midrule
\endhead
$(0,0);C_2$ & 1&1&1\\
$(0,12);R$ & 1&1&1\\
$(0,6);-$ & 17&17&12\\
$(1,10);R$ & 3&3&3\\
$(1,4);-$ & 27&25&14\\
$(2,2);-$ & 31&26&15\\
$(2,2);C_1$ & 1&1&1\\
$(2,8);R$ & 6&6&6\\
$(3,0);-$ & 16&15&13\\
$(3,0);C_1$ & 1&1&1\\
$(3,6);R$ & 12&12&10\\
$(4,10);2R$ & 1&1&1\\
$(4,4);R$ & 14&14&9\\
$(5,2);R$ & 6&6&4\\
$(5,8);2R$ & 3&3&3\\
$(6,0);R$ & 1&1&1\\
$(6,6);2R$ & 5&5&5\\
$(7,4);2R$ & 2&2&2\\
$(8,8);3R$ & 1&1&1\\
$(9,6);3R$ & 1&1&1\\
\bottomrule
\end{longtable}
\end{center}
The row $(2,2);C_1$ is the genus-$1$ curve and four isolated fixed points
missed in Dillies' list.  The eight numerical types which occur twice are
\begin{center}
\scriptsize
\begin{tabularx}{\textwidth}{@{}llllXc@{}}
\toprule
classes in~\cite{BrandhorstHofmann} & $\Fix(\sigma)$ & $\Fix(\sigma^2)$ & $\Fix(\sigma^3)$ & lower-order models & lattice datum\\
\midrule
\texttt{0.6.1.4/5} & $(2,2)$ & $8P\sqcup5R$ & $C_2$ & $W(5,4,2,1)$; $D_2$ & global\\
\texttt{0.6.2.6/7} & $(1,4)$ & $7P\sqcup4R$ & $C_3$ & $W(5,2,2,2,1)$; $D_3$ & global\\
\texttt{0.6.3.10/11} & $(2,2)$ & $6P\sqcup3R$ & $C_5$ & $W(4,2,2,2,1,1)$; $D_5$ & global\\
\texttt{0.6.3.12/13} & $(2,2)$ & $6P\sqcup3R$ & $R\sqcup C_6$ & $W(4,2,2,2,1,1)$; $D_{6+R}$ & small local\\
\texttt{0.6.4.11/12} & $(3,0)$ & $5P\sqcup2R$ & $C_7$ & $W(2,2,2,2,2,1,1)$; $D_7$ & global\\
\texttt{0.6.4.22/23} & $(2,2)$ & $4P\sqcup R$ & $R\sqcup C_6$ & invariant quartic; $D_{6+R}$ & small local\\
\texttt{0.6.5.15/16} & $(2,2)$ & $4P\sqcup R$ & $C_8$ & invariant quartic; $D_8$ & global\\
\texttt{0.6.6.6/7} & $(1,4)$ & $3P\sqcup C_1$ & $C_9$ & invariant $(2,3)$ model; $D_9$ & small local\\
\bottomrule
\end{tabularx}
\end{center}
In the second column $(p_{(2,5)},p_{(3,4)})$ is written because none of these
eight generator loci contains a curve.
Here $W(\mu_1,\ldots,\mu_s)$ denotes \cref{eq:weierstrass-six} with root
multiplicities $\mu_i$; $D_g$ is an equivariant double plane whose branch
normalization has genus $g$, and $D_{6+R}$ has branch divisor a line plus a
smooth quintic.  These equations describe the two proper powers.  In a doubled
row their fixed loci do not distinguish the two period components.  Five pairs
are the exceptional order-$6$ small local types of Brandhorst--Hofmann and
require their global lattice type; the other three are already separated by the
small local lattice type
\cite[Theorem~1.4 and Remarks~1.5--1.6]{BrandhorstHofmann}.

\begin{proposition}\label{prop:no-go}
No collection of the numbers $\eul(Y_{6,n})$, even for all $n$, classifies the
$150$ deformation classes.  Nor can the fixed loci of all powers alone do so.
There are $102$ different Euler sequences and $142$ different triples of fixed
loci.
\end{proposition}

\begin{proof}
By \cref{eq:stringy-convolution}, every stringy Euler number depends only on
$(\varepsilon_1,\varepsilon_2,\varepsilon_3)$, and the recurrence
\cref{eq:rec6} shows that the whole sequence is fixed by finitely many initial
values.  The enumeration in \cref{thm:numerical-classification} gives the
two counts.  The doubled rows above are explicit counterexamples to injectivity.  Conversely, Proposition~\ref{prop:euler-determination} shows
that the Euler sequence determines the three fixed-locus Euler numbers and
is determined by them.
\end{proof}

\medskip
After substituting the
Riemann--Hurwitz formulas \cref{eq:RHG,eq:RHF} into $\mathcal B_6$, the
coefficients of $b$ cancel ($-2b$ against $-2q_G$), and those of $a$ cancel
($-2a$ against $-(q_1+q_2)$).  The balance alone therefore does not recover the orbit decomposition from
the three fixed loci.  Proposition~\ref{prop:balance} gives the dependence on Riemann--Hurwitz.

\subsection{Quotient genera}

For order $4$, suppose that $D$ is of type II.  Given
$(r,m,k,g,n_1,n_2,a,b)$, either $R_{4,1}$ or $R_{4,2}$ gives
$g(D/\sigma)$ by Corollary~\ref{cor:qrecover}.  The numerical data are incompatible if the two answers differ or are not
integers.  Riemann--Hurwitz then
gives the number $n_2$ of branch points.  The non-zero residual detects a
positive quotient genus in type II.  A zero residual occurs both for a
pointwise-fixed maximal curve and for a type-II curve with rational quotient;
the relations alone do not distinguish these two cases.

For order $6$, one first applies
\cref{eq:D6,eq:P6,eq:T1,eq:T2,eq:T3}.  Riemann--Hurwitz checks the quotient genera for a proposed orbit
decomposition; $\mathcal B_6=0$ is a consequence of these checks.  We use these conditions to compute Calabi--Yau invariants; they do not prove
the existence of the required $K3$ lattice.

\subsection{Hodge numbers of the examples}

Using the fixed loci of the examples, we compute the Hodge numbers without
constructing each exceptional divisor.  We obtain:
\begin{center}
\begin{tabular}{@{}lll@{}}
\toprule
$K3$ model & order & $(h^{1,1},h^{2,1})$ of the threefold \\
\midrule
$w^4=Q_4$ & $4$ & $(9,27)$ \\
$y^2=x^3-a_8(t)x$ & $4$ & $(44,8)$ \\
$y^2=x^3+f_{12}(t)$, $f_{12}$ square-free & $6$ & $(29,29)$ \\
$f_{12}=t^4(t-1)^3(t+1)^3(t-\lambda)^2$ & $6$ & $(55,1)$ \\
\bottomrule
\end{tabular}
\end{center}
Collisions of roots and changes of ramification in the branch divisor change
the Hodge numbers of these Calabi--Yau threefolds.

\subsection{Higher-dimensional Euler sequences}

Let $E_{d,n}=\eul(Y_{d,n})$.

\begin{proposition}\label{prop:euler-determination}
For every $n\geq1$,
\begin{align}
E_{4,n+3}&=9E_{4,n+2}+E_{4,n+1}-9E_{4,n},\label{eq:rec4}\\
E_{6,n+4}&=12E_{6,n+3}-19E_{6,n+2}-12E_{6,n+1}+20E_{6,n}.
\label{eq:rec6}
\end{align}
For order $6$, two Euler sequences agree if and only if the triples
$(\varepsilon_1,\varepsilon_2,\varepsilon_3)$ agree.  It is enough to compare
$E_{6,2},E_{6,3},E_{6,4}$.
\end{proposition}

\begin{proof}
Let $T_d$ be convolution by $u_d/d$ on functions on $A_d$.  For a subset $H$
of $A_d$, write $\mathbf1_H$ for its indicator function.  Then
\[
 u_4=2\mathbf1_{A_4}+2\mathbf1_{2A_4}-4\mathbf1_{\{0\}},\qquad
 u_6=\mathbf1_{A_6}+2\mathbf1_{2A_6}
       +3\mathbf1_{3A_6}-6\mathbf1_{\{0\}}.
\]
The characters of $A_d$ diagonalize $T_d$.  Summing a character over each
subgroup gives the eigenvalues $9,1,-1$ for $T_4$ and $10,2,1,-1$ for $T_6$.
Their multiplicities are $1,3,12$ and $1,3,8,24$, respectively.  Thus the
annihilating polynomials are
\[
 (z-9)(z-1)(z+1),\qquad (z-10)(z-2)(z-1)(z+1).
\]
Applying them in the convolution formula proves the recurrences.

Since $E_{6,1}=24$, the next three values determine the sequence.  Conversely,
Theorem~\ref{thm:stringy-convolution} gives
\[
 \begin{pmatrix}E_{6,2}\\E_{6,3}-96\\E_{6,4}-672\end{pmatrix}
 =\begin{pmatrix}4&4&2\\64&24&8\\660&228&84\end{pmatrix}
  \begin{pmatrix}\varepsilon_1\\\varepsilon_2\\\varepsilon_3\end{pmatrix}.
\]
The matrix has determinant $-2112$, so these values recover the three
fixed-locus Euler numbers.
\end{proof}

Two values do not suffice.  Classes \texttt{0.6.3.7} and \texttt{0.6.4.11}
in the Brandhorst--Hofmann database have Euler triples $(5,-3,8)$ and
$(3,9,-12)$, respectively.  Both give
$(E_{6,2},E_{6,3})=(24,408)$, but their values of $E_{6,4}$ are $3960$ and
$3696$.  Thus equality of $E_{6,2}$ and $E_{6,3}$ does not imply equality of
the Euler sequences.

For the cyclic quartic model of Example~\ref{ex:quartic}, the first values are
\[
E_{4,1}=24,\qquad E_{4,2}=-36,\qquad E_{4,3}=-156,
\]
and \cref{eq:rec4} determines all higher dimensions.

\subsection*{Computational files}

The ancillary files contain the code and tables used in this paper.
The scripts \texttt{verify\_relations.py} and
\texttt{stringy\_euler.py} check the linear identities and compute the
convolution coefficients using exact arithmetic.  The script
\texttt{analyze\_bh\_order6.py} produces the table of $142$ numerical types,
including the three initial values which distinguish the $102$ Euler
sequences.  The file \texttt{geometric\_order6\_classification.py}
enumerates the $817$ candidates.  The necessary lattice sieve is implemented
in \texttt{order6\_discriminant\_sieve.py} using SageMath.
The scripts and tables are in \texttt{anc/computations} and
\texttt{anc/output/data}, respectively.  Instructions are given in
\texttt{anc/README.txt}.  The database input is version~2.1 of
\cite{BrandhorstHofmannData}.

\section*{Acknowledgements}
The author thanks S\l{}awomir Cynk for helpful discussions.
AI tools were used to assist with the computations and the preparation of the
manuscript.

\begingroup
\small
\enlargethispage{3\baselineskip}
\bibliographystyle{amsalpha}
\bibliography{references}

\providecommand{\bysame}{\leavevmode\hbox to3em{\hrulefill}\thinspace}
\providecommand{\MR}{\relax\ifhmode\unskip\space\fi MR }
\providecommand{\MRhref}[2]{%
  \href{http://www.ams.org/mathscinet-getitem?mr=#1}{#2}
}
\providecommand{\href}[2]{#2}
\begin{thebibliography}{DHVW86}

\bibitem[AS08]{ArtebaniSarti3}
Michela Artebani and Alessandra Sarti, \emph{Non-symplectic automorphisms of
  order 3 on {K3} surfaces}, Math. Ann. \textbf{342} (2008), no.~4, 903--921.

\bibitem[AS15]{ArtebaniSarti4}
\bysame, \emph{Symmetries of order four on {K3} surfaces}, J. Math. Soc. Japan
  \textbf{67} (2015), no.~2, 503--533.

\bibitem[Bat99]{BatyrevStringy}
Victor~V. Batyrev, \emph{Non-archimedean integrals and stringy {Euler} numbers
  of log-terminal pairs}, J. Eur. Math. Soc. \textbf{1} (1999), no.~1, 5--33.

\bibitem[BH23a]{BrandhorstHofmann}
Simon Brandhorst and Tommy Hofmann, \emph{Finite subgroups of automorphisms of
  {K3} surfaces}, Forum Math. Sigma \textbf{11} (2023), e54.

\bibitem[BH23b]{BrandhorstHofmannData}
\bysame, \emph{Finite subgroups of automorphisms of {K3} surfaces: database},
  Zenodo, version 2.1, 2023, \url{https://doi.org/10.5281/zenodo.10410278}.

\bibitem[Bor97]{Borcea}
Ciprian Borcea, \emph{{K3} surfaces with involution and mirror pairs of
  {Calabi--Yau} manifolds}, Mirror Symmetry II, AMS/IP Stud. Adv. Math.,
  vol.~1, Amer. Math. Soc., Providence, RI, 1997, pp.~717--743.

\bibitem[Bur24]{BurekBV}
Dominik Burek, \emph{Higher dimensional analogon of {Borcea--Voisin}
  {Calabi--Yau} manifolds, their {Hodge} numbers and {$L$}-functions}, Commun.
  Math. Phys. \textbf{405} (2024), no.~4, 100.

\bibitem[CdLS22]{CohnDeLaatSalmon}
Henry Cohn, David de~Laat, and Andrew Salmon, \emph{Three-point bounds for
  sphere packing}, 2022,
  \href{https://arxiv.org/abs/2206.15373}{arXiv:2206.15373}.

\bibitem[CE03]{CohnElkies}
Henry Cohn and Noam Elkies, \emph{New upper bounds on sphere packings. {I}},
  Ann. of Math. (2) \textbf{157} (2003), no.~2, 689--714.

\bibitem[CG16]{CattaneoGarbagnati}
Andrea Cattaneo and Alice Garbagnati, \emph{{Calabi--Yau} 3-folds of
  {Borcea--Voisin} type and elliptic fibrations}, Tohoku Math. J. \textbf{68}
  (2016), no.~4, 515--558.

\bibitem[CH07]{CynkHulek}
S{\l}awomir Cynk and Klaus Hulek, \emph{Higher-dimensional modular
  {Calabi--Yau} manifolds}, Canad. Math. Bull. \textbf{50} (2007), no.~4,
  486--503.

\bibitem[DHVW85]{DixonEtAl1}
Lance Dixon, Jeffrey Harvey, Cumrun Vafa, and Edward Witten, \emph{Strings on
  orbifolds}, Nuclear Phys. B \textbf{261} (1985), 678--686.

\bibitem[DHVW86]{DixonEtAl2}
\bysame, \emph{Strings on orbifolds. {II}}, Nuclear Phys. B \textbf{274}
  (1986), 285--314.

\bibitem[Dil12]{Dillies6}
Jimmy Dillies, \emph{On some order 6 non-symplectic automorphisms of elliptic
  {K3} surfaces}, Albanian J. Math. \textbf{6} (2012), no.~2, 103--114.

\bibitem[Nik80]{NikulinForms}
V.~V. Nikulin, \emph{Integral symmetric bilinear forms and some of their
  applications}, Math. USSR-Izv. \textbf{14} (1980), no.~1, 103--167.

\bibitem[Nik83]{NikulinInvolutions}
\bysame, \emph{Factor groups of groups of automorphisms of hyperbolic forms
  with respect to subgroups generated by 2-reflections. {Algebrogeometric}
  applications}, J. Soviet Math. \textbf{22} (1983), no.~4, 1401--1475.

\bibitem[Voi93]{Voisin}
Claire Voisin, \emph{Miroirs et involutions sur les surfaces {K3}},
  Ast\'erisque \textbf{218} (1993), 273--323.

\end{thebibliography}
\endgroup

\end{document}